\documentclass[11pt]{article}
\usepackage{amsmath,amsfonts,amssymb,graphicx,mathtools,flexisym, bbm, enumitem, caption, subcaption}
\usepackage{algorithm} 
\usepackage{multirow} 
\usepackage{CJK} 
\usepackage[table]{xcolor}
\usepackage{float}
\usepackage{algpseudocode}

\newcommand{\qed}{\hfill$\blacksquare$}
\newcommand{\qedsymbol}{\tag*{$\blacksquare$}}
\usepackage[margin=2cm]{geometry}
\usepackage{graphicx}
\renewcommand{\l}{\lambda}
\newcommand{\Tr}{\mathrm{Tr}}
\title{\textbf{Distributed Optimization of 
\changedSB{Bivariate}
Polynomial Graph Spectral Functions via Subgraph Optimization}}
\author{Jitian Liu, Nicolas Kozachuk, and Subhrajit Bhattacharya}
\date{}
\newtheorem{theorem}{Theorem}

\newtheorem{proposition}[theorem]{Proposition}

\newtheorem{definition}[theorem]{Definition}

\newcommand{\changedSB}[1]{#1} 

\newcommand{\todo}[1]{} 

\begin{document}
\maketitle

\begin{abstract}
We study distributed optimization of finite-degree polynomial Laplacian spectral objectives under fixed topology and a global weight budget, targeting the collective behavior of the entire spectrum rather than a few extremal eigenvalues. By re-formulating the global cost in a bilinear form, we derive local subgraph problems whose gradients approximately align with the global descent direction via an SVD-based test on the \(ZC\) matrix. This leads to an iterate-and-embed scheme over disjoint 1-hop neighborhoods that preserves feasibility by construction (positivity and budget) and scales to large geometric graphs. For objectives that depend on pairwise eigenvalue differences \(h(\lambda_i-\lambda_j)\), we obtain a quadratic upper bound in the degree vector, which motivates a ``warm-start'' by degree-regularization. The warm start uses randomized gossip to estimate global average degree, accelerating subsequent local descent while maintaining decentralization, and realizing $\sim95\%{}$ of the performance with respect to centralized optimization. We further introduce a learning-based proposer that predicts one-shot edge updates on maximal 1-hop embeddings, yielding immediate objective reductions. Together, these components form a practical, modular pipeline for spectrum-aware weight tuning that preserves constraints and applies across a broader class of whole-spectrum costs.
\end{abstract}

\vspace{-0.5em}

\vspace{-0.5em}
\section{Introduction}
%
\paragraph{Problem Statement}
Graph Laplacian matrices have long served as a pivotal role in spectral graph theory, as their eigenvalues and eigenvectors encode essential information about connectivity and stability of dynamical processes on their corresponding networks. Given a simple, connected graph $G = (V,E)$ with $|V| = n$ and weighted edges $\{w_e\}_{e\in E}$, let $A(w)$ be the adjacency matrix, $D(w) = \mathrm{diag}(A(w)\mathbbm{1})$ the degree matrix, and $L(w) = D(w) - A(w)$ the (combinatorial) Laplacian. The Laplacian spectrum $\{\lambda_i(L)\}_{i=1}^n$ (with $\lambda_1(L)\leq \lambda_2(L)\leq\cdots\leq \lambda_n(L)$) simultaneously controls, for example, spectral clustering quality, the convergence rate of consensus/diffusion (via $\lambda_2(L)$) and robustness of synchronization ($\lambda_n/\lambda_2$). This paper studies the \textit{edge-weight design} for spectral objectives under a fixed topology and a fixed total weight. Specifically, we consider the following optimization problem

\phantom{a}
\vspace{-0.1em}
\begin{equation}\label{eq:problem}
    \begin{split}
        \min_{w_1, ...,w_m} &J_G (w)\coloneqq\sum_{i, j=1 i\neq j}^n  g(\lambda_i(L(w)),\lambda_j(L(w)))\\
    &\mathrm{s.t.}\; \sum_{e\in E} w_e = W,
    \end{split}
\end{equation}

\vspace{-0.5em}\noindent
where $W$ is a positive constant and $g(\cdot, \cdot)$ is a real polynomial of two variables. The total edge weight constraint fixes that $\mathrm{tr}(L(w)) = \sum_{i}\lambda_i(L(w)) = 2\sum_e w_e = 2W$, so any improvement in $J_G$ must come from reshaping the spectrum rather than scaling it.
\changedSB{The objective function in \eqref{eq:problem} can be more compactly written as
\begin{equation} \label{eq:problem-compact}
J_G = \mathrm{Tr}\left( g(L\otimes I, I\otimes L) \right) 
\end{equation}
where $\otimes$ is the Kronecker product. A detailed justification for this representation will be given in Section~\ref{sec:problem-statement-alternative}.}

\changedSB{This type of objective functions, described as a summation over a bi-variate polynomial function evaluation of pair-wise eigenvalues of $L$, appear, for example, in context of optimization of the dispersion (variance) of the spectrum~\cite{zhang2012some} and in optimization of spectral flatness~\cite{sahin2025spectrumoptimizationdynamicnetworks}.}

\changedSB{
\paragraph{Challenges with Distributed Implementation} 
Problem~(\ref{eq:problem}) is large-scale, nonconvex, and non-local: each perturbation of $w$ moves all eigenvalues (except $\lambda_1$) through nonlinear matrix perturbation.
Although it is practically feasible to perform global, centralized optimization for small networks via numerical optimization, when the graph becomes large, the time needed for centralized optimization will be prohibitively long. Therefore, a more efficient, distributed method is desired, where optimization is performed on smaller subgraphs of $G$, which only requires the information of, and only makes changes to, the weights on the edges of the subgraph.

A (constrained) gradient descent based approach to the optimization problem is feasible.
However, evaluating each component of the gradient, $\nabla_w J_G(w)$, requires information about the entire graph.
This is unlike, for example, an objective function of the form $K_G(w) = \sum_{i} h(\lambda_i(L)) = \mathrm{Tr}(h(L))$ (with 
$h(x) = \sum_{p} c_p x^p$ 
a polynomial of degree $P$ in a single variable), the partial derivative of which with respect to the weight on the edge connecting vertices $a$ and $b$ is $\frac{K_G(w)}{\partial w_{ab}} = \mathrm{Tr}(h'(L) \frac{\partial L}{\partial w_{ab}} ) = \sum_p c_p p \mathrm{Tr}(L^{p-1} \frac{\partial L}{\partial w_{ab}})$. This can be computed only using weight information on edges within a $(P-1)$ hop neighborhood of the vertices $a$ and $b$\footnote{\changedSB{This is because the matrix $\frac{\partial L}{\partial w_{ab}}$ has non-zero entries only in locations $(a,a), (b,b), (a,b)$ and $(b,a)$, and the $a$-th (resp. $b$-th) rows/columns of $L^{p-1}$ has non-zero entries at a locations $c$ only if $c$ is at most $(p-1)$ hops away from $a$ (resp. $b$). These details are discussed further in Section~\ref{sec:acute_angle}.}}.

On the other hand, through a similar analysis, if $g(x,y) = \sum_{p,q} c_{p,q} x^p y^q$ is a polynomial of degree $P$ in $x$ and degree $Q$ in $y$, it can be verified that the gradient of the objective function $J_G$ can be computed as, 
\begin{align}
\frac{J_G(w)}{\partial w_{ab}}
&= \mathrm{Tr}\biggl( g_x (L\otimes I, I\otimes L) \left( \frac{\partial L}{\partial w_{ab}} \otimes I \right) \nonumber \\ 
&\qquad+ g_y (L\otimes I, I\otimes L) \left( I \otimes \frac{\partial L}{\partial w_{ab}} \right)  \biggr) \nonumber \\
&= \sum_{p,q} c_{p,q} \biggl( p \,\mathrm{Tr}\left( L^{p-1} \frac{\partial L}{\partial w_{ab}} \right) \, \mathrm{Tr}(L^q) \nonumber \\
&\qquad+ q \, \mathrm{Tr}(L^p) \, \mathrm{Tr}\left( L^{q-1} \frac{\partial L}{\partial w_{ab}} \right) \biggr) \label{eq:gradient-compact}
\end{align}
(this result follows from \eqref{eq:problem-compact}, and a more detailed analysis can be found in~\cite{bhattacharya2025trace}).
Because of the presence of the terms $\mathrm{Tr}(L^p)$ and $\mathrm{Tr}(L^q)$, the computation of \eqref{eq:gradient-compact} involves weights of all edges in the graph, and hence not directly amenable to distributed implementation.
}


\paragraph{Scope and Methods} This paper will primarily focus on distributed optimization of large, sparse geometric graphs, where each update touches only a small neighborhood (1-2 hops) and communicates over existing edges. The core idea is to approximate $\nabla J_G(w)$ (restricted to the local edges) with the gradient of a local problem
\begin{align*}
\min_{w_{H'}}\; &J_H(w_{H'}) \coloneqq \sum_{i,j}g(\lambda_i(L_H(w_{H'})), \lambda_j(L_H(w_{H'})))\\
&\text{s.t.}\;\sum_{e\in E_{H'}}w_e = \mathrm{const.}
\end{align*}
(where $H'$ is the subgraph whose edge are to be manipulated, and $H\supset H'$ is the subgraph where $J_H$ is supported), empirically verify the alignment between $\nabla_{w_{H'}}J_G$ and $\nabla_{w_{H'}}J_H$, solve the local problem, and embed the optimized subgraph back to the original graph. This paper also discusses a special case when $g(\lambda_i, \lambda_j) = h(\lambda_i - \lambda_j)$, where $h$ is an analytic real function. It is proved in this paper that, when $g(\lambda_i, \lambda_j) = \sum_{k\geq 0}a_k(\lambda_i - \lambda_j)^k$ with $a_{2k} \geq 0$ for any $k > 0$, there always exists some $C_0, C_h\in\mathbb{R}$ such that
\[
J_G(w) = \sum_{i,j}g(\lambda_i, \lambda_j) \leq C_0 + C_h\sum_{i,j}(d_i(w) - d_j(w))^2, 
\]
where $d(w) = A(w)\mathbbm{1}$ is the degree vector. This makes it natural to \textit{regularize} by first making vertex degrees as uniform as the topology and non-negativity allow, and then refine under the original spectral cost. Our distributed regularization method consists of two stages: 
\begin{enumerate}[noitemsep]
    \item \textit{approximation of global average degree}: adopt the well-known randomized gossip algorithm, as introduced by Boyd et al. in \cite{1638541}
    \item \textit{regularization of small neighborhoods}: manipulate the edge weights in a small subgraph such that the degree of each \textit{interior} matches the approximated average.
\end{enumerate}
Then, for the case $a_{2k} < 0$ for some $k$, the regularization may serve as a technique for \textit{warm start}: an initial guess of the edge weights that may lead to faster convergence to the minimum (achieved by centralized optimization). 

\paragraph{\changedSB{Learning-based Method for Comparison}}  In addition to optimization-based methods, a learning-based method is also developed \changedSB{for baseline comparison}. Over a training dataset generated by centralized optimization, a maximal embedding of 1-hop subgraphs is generated, and a deep neural network (DNN) model is trained to output the desired edge weights from the topology of the input 1-hop subgraph. \changedSB{While the learning-based method requires significant, conditioned training, it does not require nested local subgraphs of the form $H' \subset H$ for distributed update of the weights. However, the subgraph optimization-based approach is shown to outperform the learning-based method.}

\paragraph{Contributions} \begin{itemize}[noitemsep]
\item We formulate edge-weight optimization under fixed topology/budget for a \textit{polynomial/analytic} spectral objectives $J_G$ and develop a fully-distributed scheme that updates 1–2 hop subgraphs based on a local version of $J_G$ (namely, $J_H$) while preserving feasibility.
\item For analytic $g$ with nonnegative even-power coefficients, we derive an explicit quadratic upper bound in terms of vertex degrees, which in turn motivates a degree-regularizing scheme and serves as a warm-start surrogate when certain coefficients are negative.
\item \changedSB{We also develop a distributed learning-based approach for updating the weights locally on subgraphs, which we use as a baseline for comparison with the proposed subgraph-based optimization.}
\end{itemize}

\paragraph{Organization} The remainder of this paper is organized as follows: section 2 reviews related works on the optimization of graph Laplacian spectra are discussed. In section 3, necessary mathematical backgrounds will be presented. The problem that this paper solves will be strictly formulated in section 4, and the algorithm, along with its underlying theoretical arguments, will be presented in section 5. The special case will be discussed in section 6. Section 7 will introduce the learning-based method. The numerical experiment results for the algorithms are presented in section 8, whereas section 9 concludes the whole paper. 
\section{Related Work}
A significant body of work on Laplacian spectral optimization is founded on the principle of convex optimization. Boyd et al. established in~\cite{boyd2006convex} that many such design problems are inherently convex. The core theoretical result is that a symmetric convex function of the Laplacian's positive eigenvalues is a convex function of the graph's edge weights. This insight allows for efficient, centralized computation of optimal weights for various objectives. For instance, problems such as maximizing the algebraic connectivity ($\lambda_2(L)$) for network robustness, minimizing the spectral radius ($\lambda_n(L)$) to tolerate time delays, or finding the fastest mixing Markov chain on a graph can be formulated as semidefinite programs (SDPs) and solved to global optimality. Similarly, Shafi et al. demonstrated in~\cite{shafi2011graph} how convex optimization can be used to adjust both node and edge weights simultaneously to satisfy specific spectral bounds. In addition, a myriad of works have also emerged to maximize the algebraic connectivity of certain networks by tuning edge weights (\cite{10605101},~\cite{Tavasoli_Shakeri_Ardjmand_Rahman_2024}).

While powerful, the centralized approach faces scalability challenges, motivating distributed methods for large-scale and dynamic networks. One prominent distributed strategy involves modifying the network topology through local actions, as done in~\cite{preciado2010distributed}, where Preciado et al. proposed a distributed algorithm such that the autonomous agents add or remove edges based on limited, local information to collectively shape the global spectrum.  A key innovation of this approach is to optimize the Laplacian's spectral moments (average of Laplacian eigenvalue powers) rather than the eigenvalues directly. This is advantageous because low-order moments can be estimated using only local neighborhood data, translating a global objective into a series of localized computations. This moment-based control strategy is also central to~\cite{preciado2011spectral}, which, as a continuation of~\cite{preciado2010distributed}, further develops methods for agents to cooperatively tune their network's spectral properties without a central coordinator. 

Another class of distributed frameworks focuses on optimizing the convergence rate of discrete-time consensus protocols, which is governed by the finite condition number (namely, $\lambda_n(L)/\lambda_2(L)$). While some centralized methods have addressed this by reformulating the problem with LMIs, these are not practical for large, distributed networks. In~\cite{xu2025distributedoptimizationfinitecondition}, Xu et al. introduced a fully distributed algorithm to minimize this ratio by regulating node weights, which is a scalar quantity that is uniformly applied to all incoming information from neighboring vertices. In this context, the weighted Laplacian is defined as $L_w = \mathrm{diag}(w)L$, where $w$ is the node-weight vector, and $L$ the unweighted Laplacian. This node-weighted approach offers a significant advantage in networks where the number of nodes is much smaller than the number of edges, reducing the number of variables to optimize. However, this method inevitably breaks the symmetry of undirected networks.

Despite these progresses, the work in optimizing more general forms of Laplacian spectral functions is scant. A rational cost function of Laplacian eigenvalues is derived and studied in~\cite{sahin2025spectrumoptimizationdynamicnetworks}, yet the optimization algorithm adopted is centralized, and the problem size is relatively small. This paper, as a continuation, will establish and discuss a few methods that optimize the pairwise polynomial spectral functions in a fully decentralized manner. 
\section{Preliminaries}
\subsection{Adjacency, Degree, and Laplacian Matrices of Graphs}
Suppose one has a graph $G=(V,E)$, where $V$ is the set of vertices, and $E$ the set of edges. Then, letting $n = |V|$, $G$'s \textit{\changedSB{weighted} adjacency matrix} $A
\changedSB{\in \mathbb{R}^{n\times n}}
$ is defined via
\[
A_{ij} = \left\{ \begin{array}{ll}
    \changedSB{w_{ij}} &\text{there is an edge from vertex $i$ to vertex $j$}\\
    0 &\text{\changedSB{otherwise.}}
\end{array} \right.
\]
when $A_{ii} = 0$ for all $1\leq i\leq n$ (i.e., no self-loops), $G$ is called a \textit{simple graph}. If $A_{ij} = A_{ji}$ for any $1\leq i,j\leq n$, $G$ is called an \textit{undirected graph}. With the adjacency matrix, one may further define the \textit{degree matrix} $D$ via
\[
D \coloneqq {\rm diag}(A\mathbbm{1}_n),
\]
i.e., $D$ is a diagonal matrix whose $i$-th diagonal entry is the \textit{degree} of vertex $i$ (that is, the number of vertices that vertex $i$ is connected to). Now follows the \textit{Laplacian matrix} (which is denoted by $L$) of $G$:
\[
L\coloneqq D - A.
\]
Let $\l_1\geq \l_2\geq\cdots\geq\l_n$ be eigenvalues of $L$. Then $\l_{n-1}$ is called the \textit{algebraic connectivity} of $G$. $G$ is \textit{connected} if and only if $\l_{n-1} > 0$. All graphs appearing in this paper are assumed to be simple, undirected, and connected. Therefore, any adjacency and Laplacian matrix in this paper will be Hermitian.

\subsection{Powers of Graph Laplacian Matrices}
It can be seen from the definition of $L$ that, for $i\neq j$, $L_{ij}\neq 0$ as long as there is an edge (or a path of length 1) between vertices $i$ and $j$. Now take some positive integer $k$ and consider the matrix $L^k$. Rules of matrix multiplication tells that
\begin{align*}
    (L^k)_{ij} &= \sum_{v_1=1}^n\sum_{v_2=1}^n\cdots\sum_{v_{k-1}=1}^nL_{iv_1}L_{v_1v_2}\cdots L_{v_{k-1}j} \\
    &=\sum_{v_1=1}^n\sum_{v_2=1}^n\cdots\sum_{v_{k-1}=1}^n\prod_{\ell=0}^{k-1}L_{v_\ell v_{\ell+1}}.
\end{align*}
For the summand
$
\tilde w_{ij} = \prod_{\ell=0}^{k-1}L_{v_{\ell}v_{\ell+1}}$, \changedSB{with $v_0 = i$ and $v_k = j$},
one has that, given any sequence of vertices $(i, v_1, ..., v_{k-1}, j)$, $\tilde w_{ij}=0$ if and only if this sequence is not a path from $i$ to $j$. Therefore, in essence, the $(i,j)$-th entry of $L^k$ is the sum-of-product of all edge weights in a length $k$ path joining vertices $i$ and $j$. In particular, if the graph is unweighted, then $(L^k)_{ij}$ records the number of length $k$ paths joining $i$ and $j$.
\subsection{Problem Statement and Its Alternative Form} \label{sec:problem-statement-alternative}
Suppose that one has $n\in\mathbb{N}^+$ devices which altogether form a connected network and let a simple graph $G = (V, E, w)$ (where $w = \mathbbm{1}_{|E|\times 1}$ are the original edge weights) represent the network. Let $H\subset G$ be a induced subgraph. It is desired to manipulate edge weights in $H$ such that
\[
J_{\changedSB{G}}(w) \coloneqq \sum_{i, j \leqslant |V|, i\neq j} g(\l_i, \l_j),
\]
(where, as discussed in section 1, $g$ is a polynomial of $\lambda_i$ and $\lambda_j$, and $\l_i$ is the $i$-th smallest eigenvalue $L$) may be reduced. Let $d$ be the degree of $g$. Then one may put
\[
g(\lambda_i, \lambda_j) = \sum_{p,q=0}^{d}c_{pq}\lambda_i^p\lambda_j^q
\]
and hence
\begin{align*}
    J_G(w) &= \sum_{i,j}\sum_{p,q}c_{pq}\lambda_i^p\lambda_j^q\\ 
    &= \sum_{p,q}c_{pq}\left( \sum_i \lambda_i^p\right)\left(\sum_j \lambda_j^q \right)
    ~= \sum_{p,q}c_{pq}\mathrm{Tr}(L^p)\mathrm{Tr}(L^q)
\end{align*}
\changedSB{The right hand side of the above can also be written as $\mathrm{Tr}\left( \sum_{p,q}c_{pq}L^p \otimes L^q \right) = \mathrm{Tr}\left( \sum_{p,q}c_{pq} (L^p \otimes I) (I \otimes L^q) \right)$, giving the compact form of the objective function in terms of the Kroneker product as described in \eqref{eq:problem-compact}.}
This alternative form of $J_G(w)$ will provide convenience for theoretical analysis in the next section.
\section{Alignment between Local and Global \changedSB{Objective} Function Gradients}

\changedSB{In our distributed optimization algorithm, given an induced subgraph $H=(V_H, E_H)\subseteq G$, we want to use only the information about the weights on $E_H$ to update the weights (or a subset of the weights) on $E_H$, so as to reduce the value of the global objective, $J_G$.
The idea behind the main technical contribution of this section is to 
identify a subgraph, $H'=(V_{H'}, E_{H'})\subseteq H$, such that the gradient of an objective, $J_H$ (which depends on the weights on $E_H$) with respect to the weights on $E_{H'}$, is a subgradient of $J_G$.
}

\subsection{$k$-Hop Core of a Subgraph}
Let $H = (V_H, E_H)$ be an induced subgraph of $G$. Then one may define the \textit{$k$-hop core} of $H$ as follows:
\begin{definition}
    Let $H^c$ be the complement of $H$ in $G$. Then the \textbf{$k$-hop core} of $H$, $H'_k = (V'_k, E'_k)$, consists of the vertices that are $k$-hops away from the vertices in $H^c$, i.e.,
    $
    V'_k = \{v\in V(H): \mathrm{dist}_u(v, H^c) \geq k\}
    $
    and
    $
    E'_k = \{(u, v)\in E(H): u ,v \in V'_k(H)\},
    $
    where, 
    $
    \mathrm{dist}_u(v, H^c) = \min_{w\in H^c}\;\#\text{of hops between }v \text{ and }w. 
    $
\end{definition}

Fig.~\ref{fig:core_graph} illustrates the relationship between a 2-dimensional geometric graph $G$, its induced subgraph $H$, and the core of $H$, $H'$:
\begin{figure}[h]
    \centering
    \includegraphics[width=0.8\columnwidth, trim=400 200 400 200, clip]{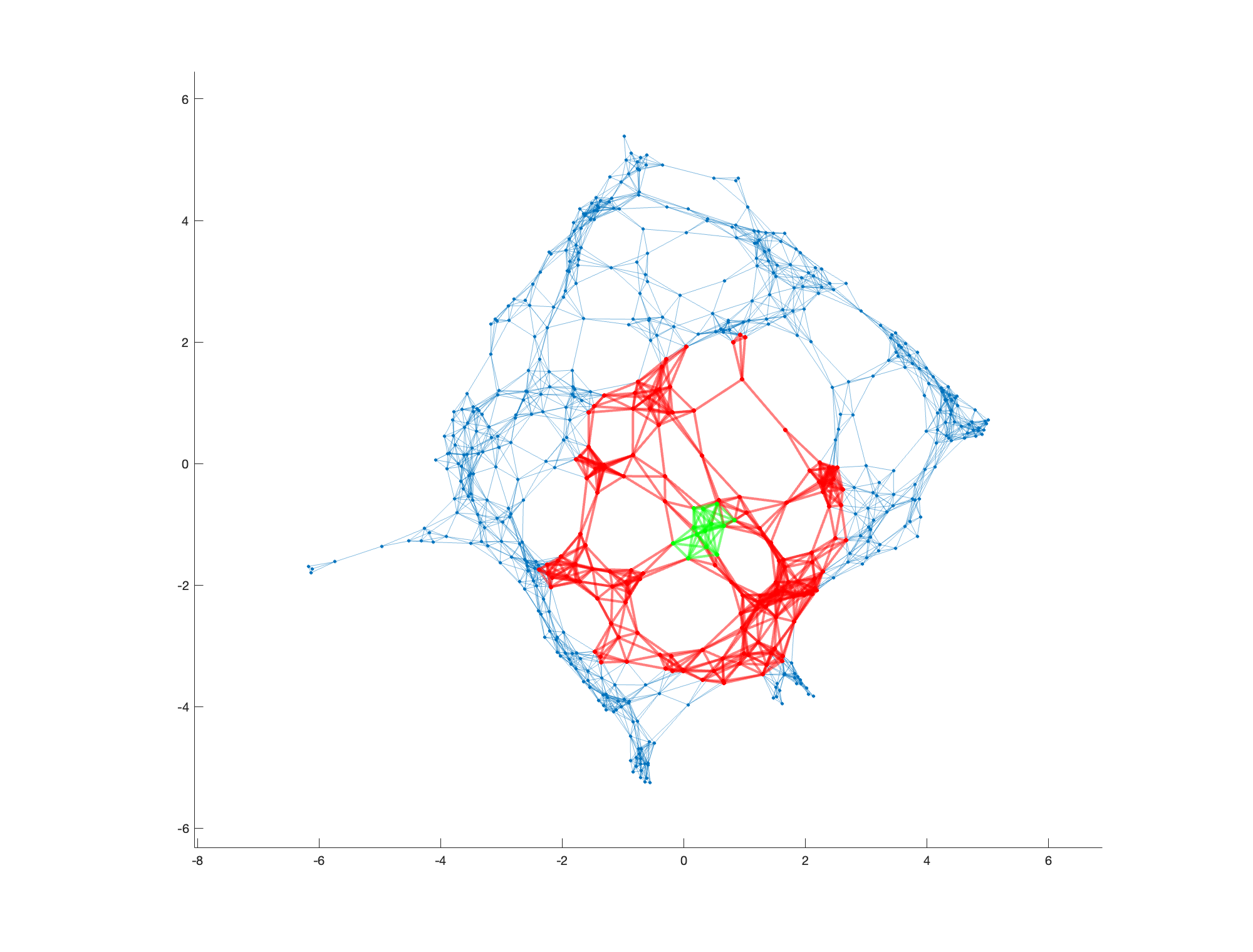}
    \caption{Visual illustration of a geometric graph \changedSB{$G$}, an induced subgraph of the graph \changedSB{$H$ (highlighted in red)}, and the subgraph's 4-hop core \changedSB{$H'$ (highlighted in green)}.
    If the degree of $g(\lambda_i, \lambda_j)$ is 4, then the weights of green edges will be tuned to minimize the local cost function supported on the subgraph consisting of the green and red edges.}
    \label{fig:core_graph}
\end{figure}

\subsection{Necessary Condition for Acute Angle between Local and Global Gradients}\label{sec:acute_angle}
Take $H$ as defined in section 5.1. To decrease $J_{\changedSB{G}}(w)$ by manipulating the weights of $E_H$, one needs to find a set of edges $E_{H'}\subseteq E_H$ (whose weights are denoted by $w_{\changedSB{{H'}}}$) such that $\nabla_{w_{\changedSB{{H'}}}} J_H$ points in the direction in which $J_G$ increases, that is, $(\nabla_{w_{\changedSB{{H'}}}}J_H)^\top (\nabla_{w_{\changedSB{{H'}}}}J_G) > 0$. Assume without loss of generality that $H$ contains the first $m < n$ vertices in $G$. Then one may \changedSB{relate the Laplcian matrices of $G$ and $H$ as}
\[
L_G = \begin{bmatrix}
L_H + D &-A\\ -A & L_{H^c} + D_c
\end{bmatrix},
\]
where $D$ is the diagonal degree matrix between $H$ and $G\setminus H$. Now consider an edge $(a,b)\in E_H$. Let $w_{ab}$ be the edge weight and 
\[
S_{ab}^2 = \frac{\partial L_H}{\partial w_{ab}}.
\]
By the definition of graph Laplacians, one knows that $S_{ab}^2$ has only 4 nonzero entries: $[S_{ab}^2]_{ab} = [S_{ab}^2]_{ba} = -1$, and $[S_{ab}^2]_{aa} = [S_{ab}^2]_{bb} = 1$ (and is therefore a rank-1 matrix). Similarly, for $L_G$,
\[
\frac{\partial L_G}{\partial w_{ab}} = \begin{bmatrix} S_{ab}^2 &0\\ 0 &0 \end{bmatrix},
\]
which may be denoted as $ \mathcal{S}_{ab}^2$. Now let $H'$ be the $d$-hop core of $H$. With the definitions of $S_{ab}^2$, the following \changedSB{condition for a positive inner product between $\nabla_{w_{\changedSB{{H'}}}}J_H$ and $\nabla_{w_{\changedSB{{H'}}}}J_G$} is discovered:
\begin{proposition}\label{prop:small_angle}
    Let
    \begin{align*}
    v_G &= \begin{bmatrix}
    \mathrm{Tr}(L_G) &
    \mathrm{Tr}(L_G^2) &
    \cdots &\mathrm{Tr} (L_G^d)
    \end{bmatrix}^\top,\\ v_H &= \begin{bmatrix}
    \mathrm{Tr}(L_H) &
    \mathrm{Tr}(L_H^2) &
    \cdots &\mathrm{Tr} (L_H^d)
    \end{bmatrix}^\top,
    \end{align*}
    and let
    \begin{align*}
    &\quad\quad\quad C = \begin{bmatrix}
    c_{11}+c_{11} &c_{12}+c_{21} &\cdots\\
    c_{21}+c_{12} &c_{22}+c_{22} &\cdots\\
    \vdots &\vdots &\ddots
    \end{bmatrix}\\ 
    z_{ab} &= \begin{bmatrix}
    \mathrm{Tr}(S_{ab}^2L_H^0)&
    2\mathrm{Tr}(S_{ab}^2L_H)&
    3\mathrm{Tr}(S_{ab}^2L_H^2)&\cdots
    \end{bmatrix}^\top.
    \end{align*}
    Taking 
    \[
    Z = \begin{bmatrix}
z_{a_1b_1}^\top\\ z_{a_1b_2}^\top\\ \vdots
\end{bmatrix},\quad (a_i, b_j)\in E(H'),
    \]
    one has that $(\nabla_{w_{\changedSB{{H'}}}}J_H)^\top(\nabla_{w_{\changedSB{{H'}}}}J_G) > 0$ only if $v_H^\top (CZ^\top ZC)v_G > 0$.
\end{proposition}
One may refer to the proof in the appendix. By consulting the proof, one knows that
\[
\nabla_{w_{\changedSB{{H'}}}}J_H = ZCv_H,\quad\nabla_{w_{\changedSB{{H'}}}}J_G=  ZCv_G.
\]
If, in particular, $ZC$'s rank is very close to 1 (i.e., for its singular value decomposition $ZC = U\Sigma V^\top$, $\sigma_1\coloneqq\Sigma_{11}$ (the largest singular value) is far larger than the remaining singular values), one may approximate the global and local gradients by
\[
ZCv_H\approx (\sigma_1v_1^\top v_H)u_1,\quad ZCv_G\approx (\sigma_1v_1^\top v_G)u_1
\]
(where $u_1$ and $v_1$ are $ZC$'s left and right singular vectors that correspond to $\sigma_1$), which strongly suggests the near-collinearlity between local and global gradients. As for the signs of $v_1^\top v_H$ and $v_1^\top v_G$, the structure of $z_{ab}$ tells that the magnitude of $Z$'s columns grow exponentially from left to right. Therefore, letting $(ZC)_j$ be the $j$-th column of $ZC$, if there is some $1\leq k\leq d$ such that
\[
|u_1^\top(ZC)_k| \gg |u_1^\top(ZC)_j|\quad\forall\;j\neq k.
\]
Then, from the equality $u_1^\top(ZC)_k = \sigma_1 (v_1)_k$, one knows that the inequality above would force $v_1 \approx e_k$ and hence $v_1^\top v_H$ and $v_1^\top v_G$ will very likely be both positive, for both $v_H$ and $v_G$ are positive.
\changedSB{These observations suggests that the condition necessary for a positive inner product between $\nabla_{w_{\changedSB{{H'}}}}J_H$ and $\nabla_{w_{\changedSB{{H'}}}}J_G$, as outlined by the proposition above, is satisfied with a high likelihood.}

\subsection{Iterative Distributed Optimization Algorithm}
Inspired by the idea in proposition~\ref{prop:small_angle}, \changedSB{we} design an iterative, distributive algorithm (summarized as algorithm~\ref{alg:cap}) that decreases $J_G(w)$ \changedSB{through subgradient descent}. In the pseudocode, the subroutines
\begin{itemize}[noitemsep]
    \item \textsc{randDraw}$(U, m)$ randomly draws $m$ elements from a set $U$,
    \item $\textsc{checkCollision}(\mathcal{H})$ returns \textbf{true} if any two subgraphs in the collection $\mathcal{H}$ have overlapping edges, and returns \textbf{true} if $\mathcal{H} = \varnothing$,
    \item \textsc{zmtx}$(H)$ computes the $Z$ matrix as defined in proposition~\ref{prop:small_angle},
    \item \textsc{svd}$(M)$ performs singular value decomposition to a matrix $M = P\Sigma Q$, and returns $P$, $\Sigma$, and $Q$,
    \item \textsc{gradDesc}$(H, H')$ computes a feasible descent direction for the optimization problem
    \[
    \min_{w_{H'}} J_H(w_{H'}),\quad\text{s.t.}\;\sum_{e\in E(H')}w_e = \text{const.}, w_{H'} \succ 0.
    \]
\end{itemize}
\begin{algorithm}
\caption{Iterative Distributed Optimization}\label{alg:cap}
\begin{algorithmic}
\Require Graph $G = (V, E, w_0)$; Number of parallel workers $m$; Cost function coefficient matrix $C$
\Ensure Optimized graph $G^*$
\State $n\gets |V|$     \Comment{Number of Vertices}
\State $U \gets \{1,...,n\}$    \Comment{Unvisited vertices}
\State $C \gets (C_{ij})_{i,j\geq 2} + (C_{ij})_{i,j\geq 2}^\top$
\State $d \gets$ number of columns in $C$
\While{$U\neq \varnothing$}
\State $\mathcal{H}'\gets\varnothing$
\While{\textsc{checkCollision}$(\mathcal{H}')=$ \textbf{true}} \Comment{Check this except right before $U$ becomes empty}
\State $V' \gets \textsc{randDraw}(U, m)$
\State $\mathcal{H}'\gets$ 1-hop neighbor graphs centered at each $v'\in V'$
\EndWhile
\State $U\gets U\setminus V'$
\For{each $H'\in\mathcal{H}'$} \Comment{Parallelized for-loop}
\State $H\gets$ $d$-hop expansion of $H'$
\State $Z\gets \textsc{zmtx}(H')$
\State $P,\Sigma,Q\gets\textsc{svd}(ZC)$
\While{$\Sigma_{11}\gg \Sigma_{kk}$ $\forall$ $k>1$ \textbf{and} $Q_1\approx e_j$ for some $1\leq j\leq d$} \Comment{Test for alignment between}
\State \Comment{local and global gradients}
\State $w_0 \gets H'.\textsc{edgeWeights}$
\State $w^*\gets \textsc{gradDesc}(H, H')$
\State $H'.\textsc{edgeWeights}\gets w^*$
\State $P,\Sigma,Q\gets\textsc{svd}(ZC)$
\EndWhile
\EndFor
\EndWhile
\end{algorithmic}
\end{algorithm}

\todo{need to emphasize that the checking of the condition for the proposition happens in line \#x in the pseudo-code, and that the computations for the check require only local information. If the condition is not satisfied, we simply move on to a different subgraph, $H$.}
\section{Special Case: When $J_G$ is a Function of Pairwise Eigen-differences}
\subsection{Relation to QP Problems for Vertex Degrees}
A significant limitation of the method in section~\ref{sec:acute_angle} is that it only applies to the case $d < \infty$. When $J_G$ is an analytic function of pairwise distances between $L_G$'s eigenvalues, i.e.,
\[
g(\lambda_i, \lambda_j) = \sum_{k = 0}^\infty a_k(\lambda_i - \lambda_j)^k,
\]
rather than optimizing over the space of edge weights, one may optimize $J_G$ by first calculating the vertex degrees that minimize a surrogate function $\tilde{J}_G$ of $J_G$ and then solve the linear equations
\[
\sum_{e\ni v}w_e = d^*_v
\]
(where $d^*_v$ is the optimal vertex degree of $v$) to calculate the desired edge weights. The following proposition dictates how one may define such a surrogate function:
\begin{proposition}
    If $g$ is analytic, then there always exists another analytic function $\tilde g$ such that 
    \[
    \sum_{i,j}g(\lambda_i,\lambda_j) \leq \sum_{i,j}\tilde g(d_i, d_j).
    \] 
    In particular, if 

    \vspace{-1em}
    \[
    g(\lambda_i, \lambda_j) = \sum_{k=0}^{\infty} \frac{a_k}{k!}(\lambda_i-\lambda_j)^k,
    \]
    then
    \begin{align*}
    &\tilde g(d_i, d_j) = a_0 +  2\sum_{k=1}^{\infty}\frac{a_{2k}}{(2k)!}\lambda_\max(L)^{2k-2}\left[(d_i-d_j)^2+ \frac{\mathrm{tr}(D)^2}{n^2}\right].
    \end{align*}
    given that $a_{2k}\geq 0$ for any $k$.
\end{proposition}
One may immediately see that this bound is \textit{sharp}: the equality is attained when the graph is complete with uniform edge weights. In practice, one may discard the constant terms and replace $\lambda_\max(L_G)$ with $2d_\max$ to search for $d\in\mathbb{R}^n$ that minimizes the function
\[
\tilde g(d_i, d_j) = \left(\sum_{k=1}^{\infty}\frac{a_{2k}}{(2k)!}d_\max^{2k-2}\right)(d_i-d_j)^2
\]
under the constraints $\sum_i d_i = \rm const.$ and $d_\max\succ d \succ 0$. It is not hard to see that $\tilde g$ is a quadratic function of $d$. Hence, the minimum of $\tilde g$ must be attained when $d_1 = \cdots = d_n$. Therefore, to solve the surrogate problem, it suffices to manipulate the edge weights such that the entire graph becomes as regular as possible, and one may achieve this by solving the following QP problem:
\[
\min_w \Vert Bw - \bar d\mathbbm{1} \Vert_2^2\quad\text{s.t.}\; w\succ 0,\quad \sum_iw_i = \text{const.},
\]
where $B$ is a $|V|\times|E|$ binary matrix with $B_{ij}=1$ if and only if the vertex $i$ is attached to the edge $j$, and $\bar d = \sum_i w_i/n$ is the average degree of all vertices of the graph.

\subsection{Distributed Local Averaging as ``Warm Start''}

In practice, if $g(\lambda_i, \lambda_j)$ can be well approximated by a finite-degree polynomial expansion, even if $a_{2k} < 0$ for some $k$, one can still initialize the edge weights by solving the surrogate problem and subsequently solve the original problem in a distributed manner following algorithm~\ref{alg:cap}.  Similarly to Section~\ref{sec:acute_angle}, one may choose multiple 1-hop disjoint subgraphs and, for each subgraph, manipulate the edge weights so that the degrees of vertices in the subgraph become homogeneous, or almost homogeneous, and hence solve the quadratic relaxation of the original problem. The distributed averaging method to be introduced requires no centralized information beyond a local estimate of $\bar d$, which one may obtain by a few rounds of neighbor-to-neighbor averaging (gossiping).

In order that the local degrees match the global degrees, at iteration $i$, for a sampled center vertex $v_i$, one shall manipulate the edges that are:
\begin{itemize}[noitemsep]
    \item between $v_i$ and $N_1(v_i)$,
    \item within $N_1(v_i)$, and
    \item between $N_1(v_i)$ and $N_2(v_i)$,
\end{itemize}
where 
\[
N_k(v_i) = \{v\in V:\mathrm{dist}_u(v, v_i) = k\}.
\]
As for estimation of the global average degree, for vertex $i$, one may start from $s^{(0)} = d$ (the current vertex degrees) and take the update rule $s^{(r + 1)} \gets \Pi s^{(r)}$, where $\Pi$ is a symmetric, doubly stochastic matrix (i.e., all row and column sums of $\Pi$ are 1). To numerically compute $\Pi$, one may randomly select $R > 0$ edges $e_1,...,e_R$ and, for each edge $e_r = (p_r, q_r)$, define 
\[
W_{e_r} = I - \frac{1}{2}(e_{p_r} - e_{q_r})(e_{p_r} - e_{q_r})^\top
\]
and let

\vspace{-1.5em}
\[
\Pi = \prod_{r=1}^R W_{e_r}.
\]
Notice that $W_{e_r}$, as an operator, averages the $e_{p_r}$-th and $e_{q_r}$-th entries of a vector. Since $P_{e_r}$ is doubly stochastic and non-negative for each $r$, $\Pi$ is also doubly stochastic. Now one may claim the following:
\begin{proposition}
    Let $\bar d = \mathbbm{1}^\top s^{(0)}/ n$, $\Pi$ as defined above, and $z^{(k)} = s^{(k)} - \bar d\mathbbm{1}$. Then one has that
    \[
    \mathbb{E} \left\Vert z^{(k+1)} \right\Vert_2^2 \leq \left( 1 - \frac{\lambda_2(L)}{2|E|}\right)^R\left\Vert z^{(k)}\right\Vert_2^2,
    \]
    and hence $s^{(k)}\to \bar d\mathbbm{1}$ almost surely.
\end{proposition}
The proof is attached in the appendix. As reference, one may also consult section II.B in~\cite{1638541} for the proof. The proposition above ensures that the randomized gossip algorithm will yield a target degree that gradually approaches the global average. The complete algorithm is summarized in algorithm~\ref{alg:gsp}.
\begin{algorithm}
\caption{Distributed Averaging via Gossiping}\label{alg:gsp}
\begin{algorithmic}
\Require Graph $G = (V, E, w_0)$; Number of parallel workers $m$; Number of gossip pairs $R$; Maximum number of iterations $\textsc{maxIter}$
\Ensure  Graph $G^* = (V, E, w_0^*)$ whose degrees are regularized
\State $n\gets |V|$     \Comment{Number of Vertices}
\State $A\gets G.\textsc{adjacency}$ \Comment{Adjacency Matrix of $G$}
\State $s\gets A\mathbbm{1}$ \Comment{Vector of vertex degrees}
\For{$i=1,...,\textsc{maxIter}$}
\For{$r = 1,...,R$} \Comment{Gossiping}
\State $(p, q)\gets$ random edge drawn from $E$
\State $s'\gets (s[p] + s[q])/2$
\State $s[p]\gets s'$
\State $s[q]\gets s'$
\EndFor
\State $\mathcal{H}\gets\varnothing$
\While{$\textsc{checkCollision}(\mathcal{H})=\textbf{true}$}
\State $V'\gets\textsc{randDraw}(V,m)$
\State $\mathcal{H}\gets$ 2-hop neighborhoods centered at each $v'\in V'$
\EndWhile
\For{each $H\in\mathcal{H}$}\Comment{Parallelized For-Loop}
\State $v\gets$ center of $H$
\State $E'\gets$ all edges that contain $v$ and $N_1(v)$
\State $B\gets$ incidence matrix between $\{v\}\cup N_1(v)$ and $E'$
\State $w^*\gets \mathrm{arg}\min_{w(E')}\Vert Bw(E') - s[\{v'\}\cup N_1(v')] \Vert_2^2$
\State $E'.\textsc{weights}\gets w^*$
\EndFor
\EndFor
\end{algorithmic}
\end{algorithm}

\section{Learning-Based Distributed Optimization}
\label{sec:learning}
\changedSB{As a baseline method for comparing against the proposed subgraph optimization method, in}
this section \changedSB{we present} a novel learning-based decentralized approach to spectrum optimization for network graphs. We propose a DNN model trained on localized outcomes from centralized spectrum optimization. The model learns to predict edge weight adjustments based on local graph topologies \changedSB{and weights}, enabling decentralized control without global network knowledge. Our methodology involves generating diverse graph datasets, embedding maximal subgraphs to standardize inputs, and iteratively refining the DNN architecture for optimal performance. 

\subsection{Methodology}

Given an objective function $J_G(\mathbf{w}_G)$ (where $\mathbf{w}_G$ are the weights on the edges of the entire graph, $G$), the centralized approach performs a gradient descent of $J_G$. With the aim of implementing this optimization in a decentralized manner, in this section we instead propose to learn one-shot optimization of $\widetilde{J}_H(\mathbf{w}_H)$, where $\mathbf{w}_H$ are the weights on a subgraph $H\subseteq G$, such that descending the subgraph gradient also results in decrease of the original objective function. This allows us to tune the edge weights on subgraphs only, and this achieves decentralized optimization.

\begin{figure*}
\centering
\includegraphics[width=0.9\linewidth]{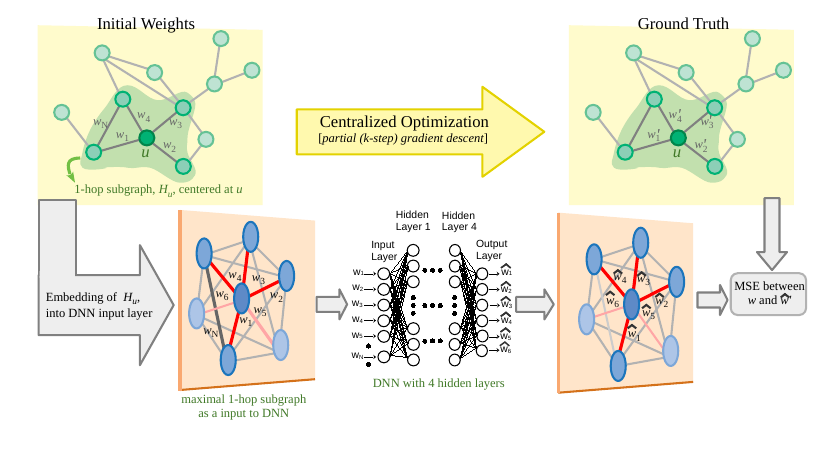}
\caption[Schematic of Overall Model Structure]{\label{fig:methodology} \changedSB{\emph{Learning-based Method:}} Schematic illustrating the model prediction pipeline, where input graph data undergoes a 1-hop maximal subgraph embedding of graph edge weights for node $u$, which is then fed into a DNN with 4 layers. The initial weights inputted ($w_1$ to $w_6$) are the edges connected to the center node(indicated by a red edge, where bright red are the present edges and the lighter red are the edges with a weight of 0). The following edges inputted are the edges within the 1-hop neighborhood that are not connected to the center node (indicated by a gray edge, where again the dark gray represents present edges and the lighter gray are the edges with a weight of 0). The model outputs predicted edge weights for the edges connected to the center node to perform spectrum optimization, and accuracy is evaluated using Mean Squared Error (MSE) by comparing predictions against the ground truth edge weights, which consist of the initial weights after undergoing centralized optimization.}
\end{figure*}

In our trained DNN model we would provide the weights on the edges of a subgraph, $\mathbf{w}_H$, as input, and the expected output of the DNN would be a set of new weights, $\mathbf{w}'_H$, on the edges of the subgraph, which is expected to reduce the global objective function, $J(\mathbf{w}_G)$. This process is demonstrated in Figure~\ref{fig:methodology}. Performing this change in weights over multiple subgraphs (each node in the graphs 1-hop/2-hop neighborhood) is expected to reduce the global objective in steps. 

\subsection{Maximal Subgraph Embedding} \label{sec:max_embed}

The DNN model's inputs consists of the 1-hop and 2-hop neighborhood subgraphs for each node within each of the generated graphs. To enable decentralized spectrum optimization, the input was intentionally restricted to this local information available to individual nodes, rather than the entire graph as is with centralized spectrum optimization.

To ensure consistent training across the varying sizes of 1-hop and 2-hop neighborhoods, the concept of a maximal subgraph embedding was introduced. This framework embeds each subgraph into a fixed-size structure, preserving the topological variations of the original subgraphs while standardizing the input dimensions. The maximal subgraph serves as a unified embedding space, generalizing all subgraphs within the dataset into a single, consistent graph structure. This facilitates efficient training and optimization by addressing the heterogeneity of subgraph sizes and enabling scalable processing across the dataset.


\begin{figure*}
\centering
\includegraphics[width=.75\linewidth]{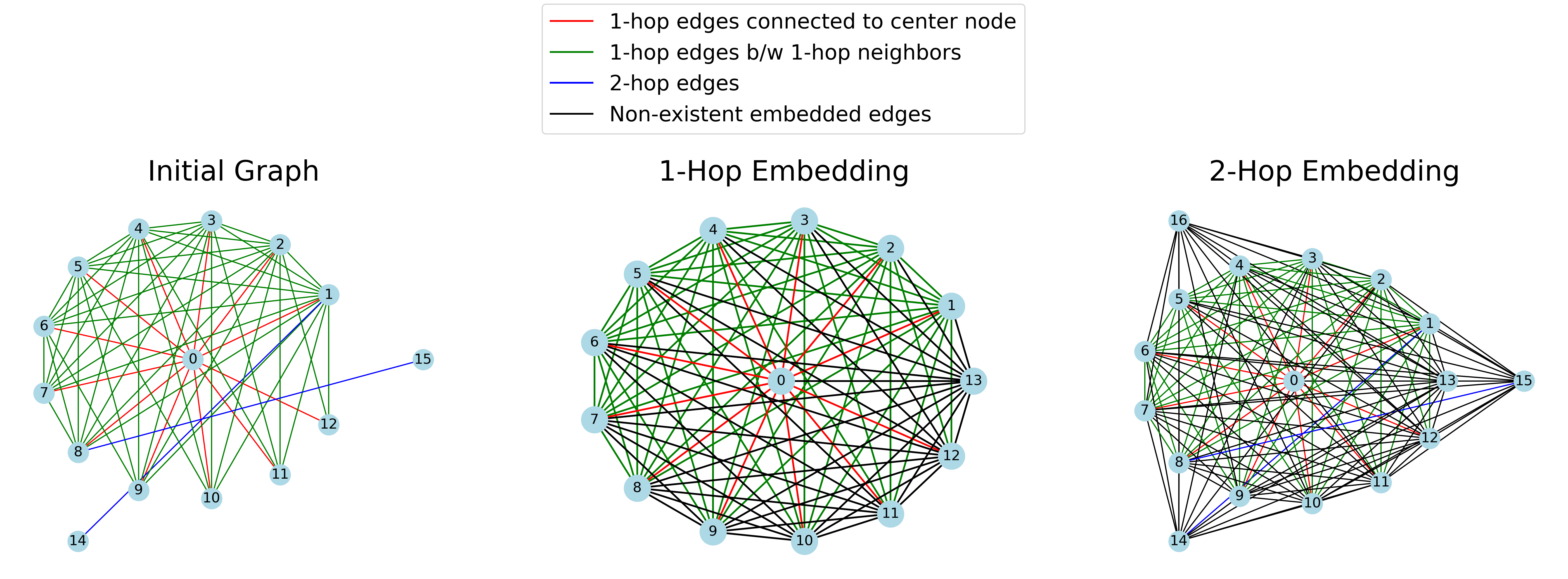}
\caption[Maximal Subgraph Embedding]{\label{fig:maximal_embedding}Example of 1-hop and 2-hop maximal embeddings of a graph with 12 1-hop subgraph nodes and 2 2-hop subgraph nodes. The left graph is the original graph. The center graph is the 1-hop maximal embedding, where the 1-hop maximal subgraph has 13 nodes. The right graph is the 2-hop maximal embedding, where the maximal subgraph has 13 1-hop nodes and 3 2-hop nodes. The red edges are 1-hop edges connected to the center node, the green edges are the 1-hop edges that are not connected to the center node, the blue edges are 2-hop edges(excluding the 1-hop edges), and the black edges are edges that do not exist in the graph but are present in the maximal embedding.}
\end{figure*}

\changedSB{For a given class of graphs, we assume that we know, a priori, the structure of the required maximal subgraph in which any 1-hop and 2-hop subgraph of the graphs can be embedded. In practice,}
the process of generating the maximal subgraph embeddings \changedSB{can be achieved} by examining each node in every graph \changedSB{of the dataset}, treating it as the central node, and extracting its 1-hop and 2-hop neighborhood subgraphs.
Each individual subgraph is embedded into the topology of this maximal subgraph, ensuring a consistent graph structure for both the input and output \changedSB{layers} of the model. This embedding process results in graphs with the same node counts and overall structure, but with the variation amongst the embeddings existing in the edge weights themselves, where a weight of 0 indicates the absence of a corresponding edge in the original subgraph when mapped to the maximal subgraph, preserving the structure of the initial subgraph. This approach maintains structural uniformity while preserving the unique characteristics of each subgraph through edge weight differences. By looking at Figure~\ref{fig:maximal_embedding}, one can see how a 14 node 2-hop neighborhood subgraph, with 12 1-hop nodes and 2 2-hop nodes, is embedded into a 12 node 1-hop subgraph, as well as embedded in a dual 13 node 1-hop and 3 node strictly 2-hop neighborhood subgraph.

\subsection{Model Generation and Structure }
Given the demonstrated learning capacity and architectural flexibility of deep neural networks, it was decided that a DNN based approach for this task could provide the best results. 
The model consisted of an input vector that is the size of the maximum number of one hop neighbors, a value for each present edge and zeros for non-existing edges, and a output of the same size representing the edge weight modifications. The model consists of four hidden layers with dimensions 512, 256, 128 and 64, and corresponding activation layers of Swish, ReLU, ReLU, and Swish (see Figure~\ref{fig:methodology}).

\subsection{Applying Predictions for Optimization}
\label{sec:apply_opt}

\changedSB{The model was trained on a set of training graphs by}
stochastically iterating over each node in a graph, one at a time, \changedSB{and updating the model parameters through backpropagation to minimize the mean squared error between the predicted weights and weights computed by a few steps of the centralized optimization}.
For each node, the current connected edge weights of its 1-hop neighborhood are used as inputs to the model, and the output provides the predicted spectrum optimized edge weights. To ensure consistency with the problem constraints, specifically that the minimum edge weight must be 0.1 and that the total sum of edge weights before and after optimization must remain unchanged, the outputted edge weights are normalized to satisfy these constraints before applying the edge weights manipulations to the network graph. This normalization process preserved the structural integrity of the graph while adhering to the imposed constraints.

After \changedSB{each training epoch using a batch of training graphs, the model was tested on a new batch of graphs by} applying the predicted optimized edge weights 
\changedSB{to the 1-hop neighborhood of each node of the test graphs.}
\changedSB{It was observed that with each subsequent training epoch, the model performed better at updating the weights on the 1-hop edges, for achieving a global objective value that is closer to the one achieved by centralized optimization (Figure~\ref{fig:dopr_vs_size}. Thus the} outcome demonstrates the model's ability to effectively perform decentralized spectrum optimization, achieving the desired improvements in network resilience while adhering to the necessary constraints.
\section{Numerical Experiment, Results, and Observations}
The cold-start, warm-start, and learning-based methods were tested by solving the problem
\begin{align*}
\min_{w} J_G(w) &= \sum_{i,j=1}^n(\lambda_i - \lambda_j)^4 - (\lambda_i - \lambda_j)^2,\\ &\text{s.t.}\;\sum_{e\in E}w_e = \text{const.},\; w\succ 0.
\end{align*}
According to section 4, the cost function may be re-written in the bilinear form
\[
J_G(w) = \sum_{p,q=1}^5 \mathbf{C}_{pq}\mathrm{Tr}(L^{p-1})\mathrm{Tr}(L^{q-1}),
\]
where
\[
\mathbf{C} = \begin{bmatrix}
0 &0 &-1 &0  &1\\
0 &2 & 0 &-4 &0\\
-1 &0 &6 &0 &0\\
0 &-4&0 &0 &0\\
1 & 0 &0 &0 &0
\end{bmatrix}.
\]
To quantify the performance of our methods, we introduce an evaluation metric named the \textit{Decentralized Optimization Performance Ratio} (DOPR). The DOPR is defined as the ratio of the objective value reduction achieved by distributed optimization to that of centralized optimization, expressed as:
\[
\mathrm{DOPR} = \frac{J_0 - J_d}{J_0 - J^*},
\]
where $J^*$ is the value of the cost function after optimization, $J_0$ the value associated with uniform edge weights, and $J_d$ the value after distributed optimization. In the experiment, 350 random geometric graphs with 600 nodes are generated with uniform edge weights. To generate each graph, 600 points are sampled on $(0,1)^2$ based on the uniform distribution, and any pair of points $(p, q)$ with $\Vert p-q\Vert_2 < 0.08$ is assigned an edge. As a baseline, centralized optimization is run for each graph on MATLAB via the black-box solver \texttt{fmincon()} to directly solve (\ref{eq:problem}).

\subsection{Cold-Start vs Warm-Start}
\begin{figure}
    \centering
    \includegraphics[width=\linewidth]{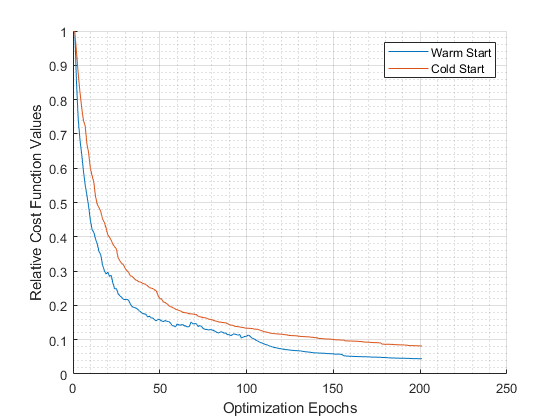}
    \caption{Descent curves of warm start versus cold start on the same graph. The $x$-axis marks the optimization epochs, whereas the $y$-axis records the percentage of decrease with respect to centralized method (namely, $1 - \rm DOPR$). The value $1$ corresponds to the fact that the objective value does not decrease at all, where as the value $0$ indicates that the objective value reaches the level after centralized optimization.}
    \label{fig:cold_vs_warm}
\end{figure}
\begin{figure}
    \centering
    \includegraphics[width=\linewidth]{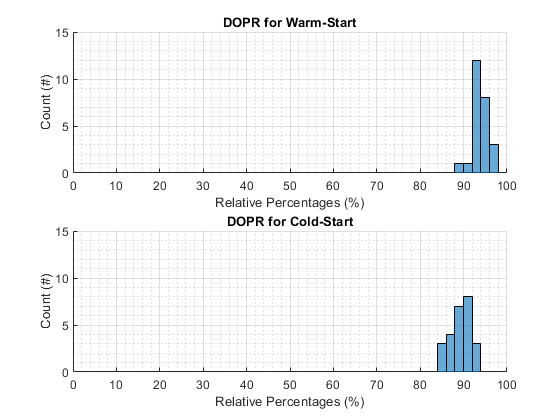}
    \caption{Histograms for DOPRs of cold start and warm start methods. It can be noticed from the figure that the DOPR of warm-start may reach $\sim95\%{}$ in average, whereas cold-start may only reach $\sim90\%{}$.}
    \label{fig:hist_warm_vs_cold}
\end{figure}

Since $J_G$ is a function of pairwise eigen-difference, the performances of subgraph optimization with and without local degree averaging are compared across the same graph. Both methods have been ran for 200 iterations. For the cold start method, in all 200 iterations, 8 disjoint subgraphs are randomly selected from the entire graph and the edge weights are tuned to solve local problems.For the warm start method, the first 100 iterations are run to locally average the vertex degrees, and the last 100 iterations are for distributed optimization. One may refer to Fig.~\ref{fig:cold_vs_warm} for a visual comparison between descent curves for both methods. Also, the cold start and warm start methods are compared across 25 graphs. One may refer to Fig.~\ref{fig:hist_warm_vs_cold} for a comparison between their performances with respect to the centralized approach. 
\begin{figure}
    \centering
    \includegraphics[width=\linewidth]{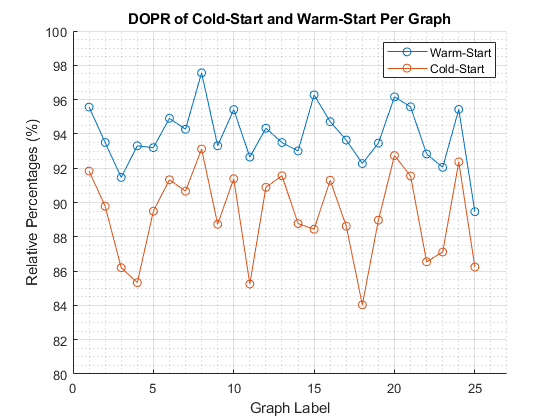}
    \caption{Comparison between DOPRs of cold-start and warm-start methods for each graph in the testing dataset. One may see from the curves that, in our test case, the warm-start method is \textit{strictly superior} to the cold-start method.}
    \label{fig:placeholder}
\end{figure}

\subsection{Learning-Based vs Cold-Start}
\begin{figure}
    \centering
    \includegraphics[width=\linewidth]{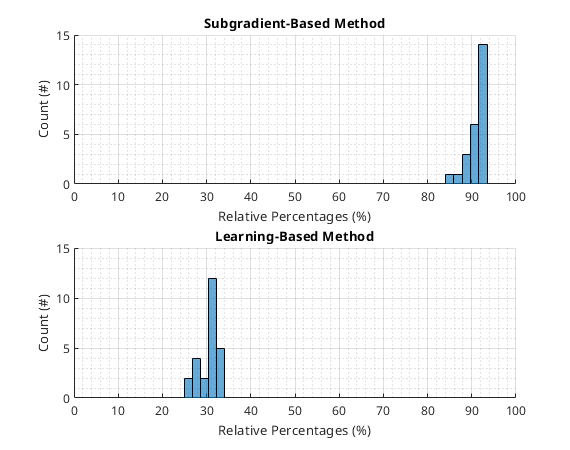}
    \caption{The histograms of DOPRs for the cold start and learning-based methods. As indicated in the figure, the top histogram plots the DOPRs of cold start method, and the bottom the learning-based method.}
    \label{fig:subgrad_vs_learning}
\end{figure}

The majority of the 300 generated graphs serves as training dataset to train the DNN proposed in section~\ref{sec:learning}. Among the 300 graphs, 275 are used for training and 25 for testing. For comparison, the cold-start distributed optimization method is also run on the test set, and the DOPRs for both methods are compared, as shown in Fig.~\ref{fig:subgrad_vs_learning}. It is shown that, the average DOPR of learning-based methods will reach around 30\%{}, whereas, in contrast, the cold start method can reach around 90\%{} in average. 

\subsection{More Observations}
\begin{figure}
    \centering
    \includegraphics[width=\linewidth]{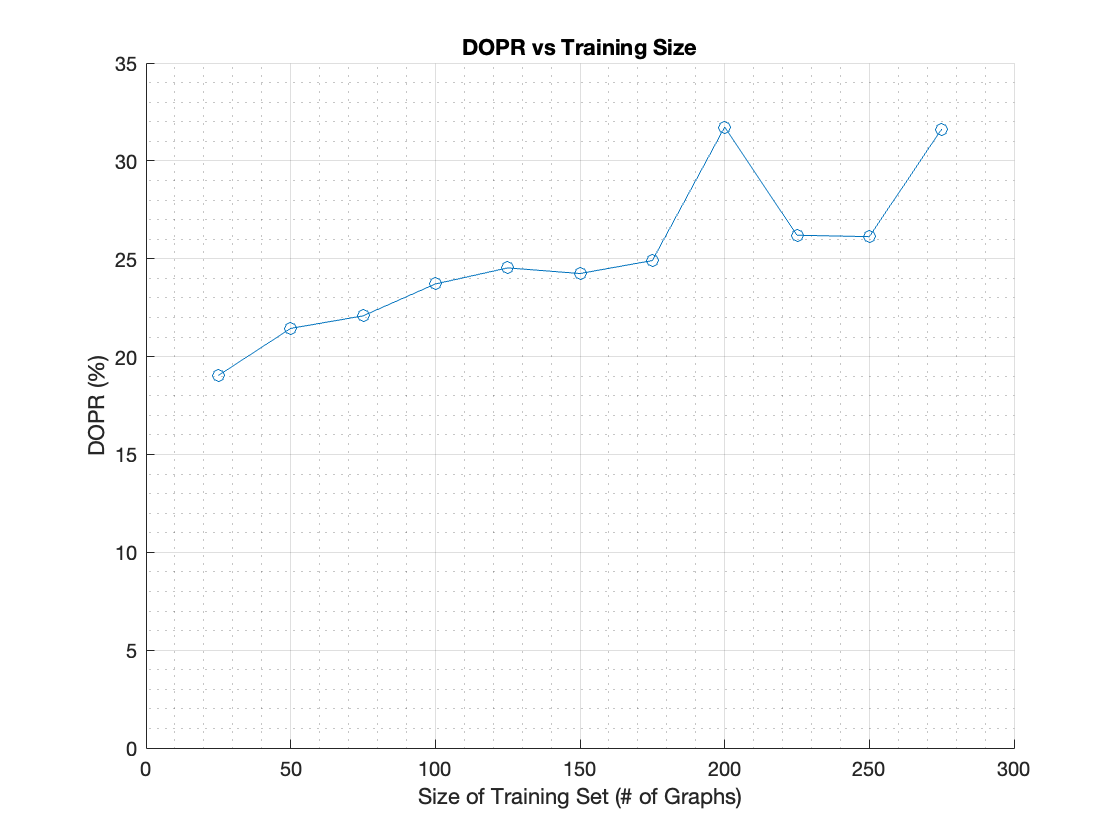}
    \caption{DOPR of the trained model vs the size of training data.}
    \label{fig:dopr_vs_size}
\end{figure}

Using this DOPR metric and performing both centralized and decentralized optimization, it was calculated that the DOPR value for decentralized optimization is 0.29. This means that decentralized optimization achieves an objective function value decrease that is 29\% that of what the centralized optimization can achieve. Moreover, training the network with dataset of different sizes, an increase in the DOPR was observed, as shown in Fig.~\ref{fig:dopr_vs_size}.

\section{Conclusion}

In this paper, we present a few distributed frameworks that optimize finite-degree polynomial Laplacian spectral functions, rather than a particular eigenvalue (like $\lambda_2$ or $\lambda_n$), under fixed topology and total weight budget constraints. By expressing the global cost function in a bilinear form, we derived local subgraph problems whose gradients (approximately) align with the global descent direction via an SVD test on the $ZC$ matrix. This led to an efficient iterate-and-embed scheme over disjoint 1-hop neighborhoods that preserves feasibility and scales to large geometric graphs. For cost objectives that are functions of pairwise eigenvalue differences, we derived a quadratic upper bound of the degree vector, suggesting an inexpensive degree-regularization ``warm start'', followed by local optimization. Finally, we developed a learning framework that predicts one-shot edge updates on maximal 1-hop embeddings and showed that it could also rather effectively reduce the objective values.

Our results extend the centralized convex formulations for Laplacian spectrum design (e.g., SDP formulations for fastest mixing, effective resistance, and algebraic connectivity) and the body of distributed methods that focus on spectral moments using only local operations. Moreover, our warm-start stage leverages randomized gossip to estimate global degree targets with provable convergence that depends on $\lambda_2$, aligning with classical consensus theory. Most importantly, our work bridges decentralized optimization and a broader class of network design problems concerning the collective behavior of the entire Laplacian spectrum (for example, the resonance-based vulnerability in networked second-order dynamics, as discussed in~\cite{sahin2025spectrumoptimizationdynamicnetworks}).

As for future work, a natural next step is to extend the class of cost functions from finite-degree polynomials to pairwise rational (and more generally, analytic) functions of the Laplacian spectrum. Also, a significant drawback of the iterate-and-embed scheme is that the radius of the subgraph (on which the local cost function is supported) significantly relies on the degree of the $g(\lambda_i, \lambda_j)$--the higher $\deg g$ is, the larger the subgraph needs to be. This will motivate further theoretical explorations about in what circumstances (form of objective function, topology of subgraphs, etc.) local optimization will lead to global improvement.

\changedSB{
\section{Acknowledgments}
We gratefully acknowledge the support of Air Force Office of Scientific Research (AFOSR) award \# FA9550-23-1-0046.
}

\bibliographystyle{plain} 
\bibliography{refs} 

\section*{Appendix}
\textbf{Proof to Proposition 2.} Now, in order to explicitly compute the partial derivative $\partial J_G/\partial w_{ab}$, it remains to derive the derivative of the map $f: M\mapsto \Tr(M^p)$. To compute the total derivative $D f$, let $\delta M$ be a small perturbation to $M$. Then
\begin{align*}
    f(M + \delta M) &= \Tr\left[(M + \delta M)^p \right]\\
    &= \Tr\left[M^p + \sum_{k=0}^{p-1} M^k(\delta M) M^{p-k-1} \right]\\ &\quad+ o(\delta M^2)\\
    &= f(M) + \sum_{k=0}^{p-1}\Tr \left[M^k(\delta M) M^{p-k-1} \right]\\ &\quad+ o(\delta M^2)\\
    &= f(M) + \sum_{k=0}^{p-1}\Tr\left(M^{p-1}\delta M\right) + o(\delta M^2)\\
    &= f(M) + p\Tr(M^{p-1}\delta M) + o(\delta M^2).
\end{align*}
Hence the derivative is $Df: \mathbb{R}^{n\times n}\to \mathbb{R}$ with $E\mapsto p\Tr(M^{p-1}E)$. Therefore, by the chain rule,
\begin{align*}
\frac{\partial f}{\partial w_{ab}} 
&= p\mathrm{Tr}\left(M^{p-1}\frac{\partial M}{\partial w_{ab}}\right),
\end{align*}
and
\begin{align*}
\frac{\partial J_G}{\partial w_{ab}} &= \sum_{p,q=1}^{d}c_{pq}\biggl[
p\mathrm{Tr}(L_G^{p-1}\mathcal{S}_{ab}^2)\mathrm{Tr}(L_G^q)\\
&\quad\quad\quad\quad\quad+
q\mathrm{Tr}(L_G^{q-1}\mathcal{S}_{ab}^2)\mathrm{Tr}(L_G^p)
\biggr]\\
&=\sum_{p,q=1}^{d}c_{pq}\biggl[
p\mathrm{Tr}(L_G^{p-1}\mathcal{S}_{ab}^2)\mathrm{Tr}(L_G^q)\\
&\quad\quad\quad\quad\quad+
q\mathrm{Tr}(L_G^{q-1}\mathcal{S}_{ab}^2)\mathrm{Tr}(L_G^p)
\biggr]\\
&=\sum_{p=1}^{d}\sum_{q=1}^{d}(c_{pq} + c_{qp})\left[
p\mathrm{Tr}(L_G^{p-1}\mathcal{S}_{ab}^2)\mathrm{Tr}(L_G^q)
\right]\\
&=\sum_{p=1}^{d} p\mathrm{Tr}(L_G^{p-1}\mathcal{S}_{ab}^2)
\sum_{q=1}^{d}(c_{qp}+c_{pq})\mathrm{Tr}(L_G^q).
\end{align*}
Now, if both vertices $a$ and $b$ are $d$-hops away from the boundary between $H$ and $G\setminus H$, then every path of length not larger than $d$ that connects $a$ and $b$ must lie entirely in $H$. Then, using the notation in section 3, one has that $\tilde w^G_{ij} = \tilde w^H_{ij} \implies (L_H^{p})_{ab} = (L_G^{p})_{ab}$ for any $0\leq p \leq d$. Thus
\[
\mathrm{Tr}(L_G^{p-1}\mathcal{S}_{ab}^2) = \mathrm{Tr}(L_H^{p-1}S_{ab}^2).
\]
Letting 
\[
C = \begin{bmatrix}
c_{11}+c_{11} &c_{12}+c_{21} &\cdots\\
c_{21}+c_{12} &c_{22}+c_{22} &\cdots\\
\vdots &\vdots &\ddots
\end{bmatrix}
\]
and 
\[
z_{ab} = \begin{bmatrix}
\mathrm{Tr}(L_H^0S_{ab}^2)&
2\mathrm{Tr}(L_HS_{ab}^2)&
3\mathrm{Tr}(L_H^2S_{ab}^2)&\cdots
\end{bmatrix}^\top,
\]
one may finally put that
\[
\frac{\partial J_G}{\partial w_{ab}} = z_{ab}^\top C v_G,
\]
where
\[
v_G = \begin{bmatrix}
\mathrm{Tr}(L_G) &
\mathrm{Tr}(L_G^2) &
\cdots &\Tr (L_G^d)
\end{bmatrix}^\top.
\]
Let
\[
Z = \begin{bmatrix}
z_{a_1b_1}^\top\\ z_{a_1b_2}^\top\\ \vdots
\end{bmatrix},\quad (a_i, b_j)\in E_H',
\]
where $E_H'$ are the edges in $H$ that are at least $d$-hops away from the boundary. One may put that
\[
\nabla_{w_{\changedSB{{H'}}}}J_G = ZCv_G
\]
and, similarly, $\nabla_{w_{\changedSB{{H'}}}}J_H = ZCv_H$. Thus, \[(\nabla_{w_{\changedSB{{H'}}}}J_H)^\top (\nabla_{w_{\changedSB{{H'}}}}J_G) = v_H^\top(C^\top Z^\top Z C)v_G.\]  This completes the proof.\qed

\bigskip
\noindent
\textbf{Proof to Proposition 3.} 
Notice first that $g(\lambda_i,\lambda_j)$ can be written as a function of a single variable $h(\lambda_i - \lambda_j)$. Also, by letting $M = D\otimes I - I\otimes D$ and $N = A\otimes I - I\otimes A$, one has that 
\[
J_G(w) = \mathrm{tr}(h(M - N)),
\]
and the problem now becomes to find a surrogate function $\tilde g$ such that $\mathrm{tr}(\tilde h(M - N)) \geq \mathrm{tr}(h(M-N))$. Since $h$ is analytic, by definition, one then has a sequence of coefficients $\{a_n\}\subset\mathbb{R}$ such that, for any Hermitian $X$,
\[
h(X) = \sum_{n=0}^{\infty}\frac{a_n}{n!}X^n \implies \mathrm{tr}(h(X)) = \sum_{n=0}^{\infty}\frac{a_n}{n!}\mathrm{tr}(\Lambda(X)^n),
\]
where $\Lambda(X)$ is the diagonal matrix containing $X$'s eigenvalues. Now, if $X = M-N$, since the spectrum of $M-N$ is symmetric about zero, whenever $n$ is odd, one will have that $\mathrm{tr}(\Lambda(M-N)^n) = 0$. Thus all one needs to focus on are the even summands, i.e.,
\begin{align*}
    \mathrm{tr}(h(M-N)) = a_0n^2 + \sum_{k=1}^{\infty}\frac{a_{2k}}{(2k)!}\mathrm{tr}(\Lambda(M-N)^{2k}).
\end{align*}
Now, for each $k$,
\begin{align*}
    \mathrm{tr}(\Lambda(M-N)^{2k}) = \Vert (M-N)^k \Vert_F^2,
\end{align*}
which can be upper bounded as follows:
\begin{align*}
    \Vert (M-N)^k \Vert_F &\leq \Vert M-N\Vert_2\Vert (M-N)^{k-1} \Vert_F\\ &\leq\cdots\\ &\leq \Vert M-N\Vert_2^{k-1}\Vert M-N \Vert_F.
\end{align*}
Squaring on both sides, one has that
\begin{align*}
    \mathrm{tr}(\Lambda(M-N)^{2k}) &\leq \Vert M-N\Vert_2^{2k-2}\mathrm{tr}[(M-N)^2]\\ &= \lambda_\max(L)^{2k-2}\mathrm{tr}((M-N)^2).
\end{align*}
Notice also that
\begin{align*}
\mathrm{tr}((M-N)^2) &=
    \sum_{i=1}^{n^2}(\lambda_i(M-N))^2 \\ &= \sum_{i,j}[\lambda_i(L) - \lambda_j(L)]^2\\ &= \sum_{i, j}[\lambda_i(L)^2 + \lambda_j(L)^2 - 2\lambda_i(L)\lambda_j(L)]\\ &= n\sum_i\lambda_i(L)^2 + n\sum_j\lambda_j(L)^2\\ &\quad-2\sum_{i, j}\lambda_i(L)\lambda_j(L)\\
    &= n\sum_i\lambda_i(L)^2 + n\sum_j\lambda_j(L)^2 \\
    &\quad-2\left[\sum_i\lambda_i(L) \right]\left[\sum_j\lambda_j(L) \right]\\ &= 2n\mathrm{tr}(L^2) - 2\mathrm{tr}(L)^2.
\end{align*}
Then, because $D_{ii} = \sum_j A_{ji} = \sum_jA_{ij}$, one has that $\mathrm{tr}(AD) = \mathrm{tr}(DA) = 0$. Therefore,
\begin{align*}
    \mathrm{tr}(L^2) = \mathrm{tr}[(D-A)^2] = \mathrm{tr}(D^2 + A^2),
\end{align*}
and, since all diagonal entries of $A$ are zero,
\begin{align*}
    \mathrm{tr}(L)^2 = [\mathrm{tr}(D-A)]^2 = [\mathrm{tr}(D) - \mathrm{tr}(A)]^2  = \mathrm{tr}(D)^2.
\end{align*}
Hence,
\begin{align*}
\mathrm{tr}((M-N)^2) &= 2n\mathrm{tr}(L^2) - 2\mathrm{tr}(L)^2\\ &= 2n\mathrm{tr}(D^2) - 2\mathrm{tr}(D)^2 + 2n\mathrm{tr}(A^2)\\ &= \mathrm{tr}(M^2) + \mathrm{tr}(N^2).
\end{align*}
Note then that
\begin{align*}
\mathrm{tr}(N^2) &= 2n\mathrm{tr}(A^2) = 2n\Vert A\Vert_F^2\\ &= 2n\sum_{i,j}A_{ij}^2 = 2n\sum_i\left(\sum_jA_{ij}^2\right)\\ &\leq 2n\mathrm{tr}(D^2).
\end{align*}
Thus,
\begin{align*}
    \mathrm{tr}((M-N)^2) &\leq \mathrm{tr}(M^2) + 2n\mathrm{tr}(D^2).
\end{align*}
Knowing that
\begin{align*}
\mathrm{tr}(M^2) &= \sum_{i,j}(d_i - d_j)^2\\ 
&= 2n\sum_{i}d_i^2 -2\left(\sum_{i}d_i \right)^2\\  &=2n\mathrm{tr}(D^2) - 2\mathrm{tr}(D)^2.
\end{align*}
one has that
\[
 \mathrm{tr}((M-N)^2) \leq 2\mathrm{tr}(M^2) + 2\mathrm{tr}(D)^2.
\]
To sum up,
\begin{align*}
    \mathrm{tr}(\Lambda(M-N)^{2k}) \leq 2\lambda_\max(L)^{2k-2}[\mathrm{tr}(M^2) + \mathrm{tr}(D)^2].
\end{align*}
and hence, in total, if $a_{2k} > 0$ for any $k$,
\begin{align*}
    \mathrm{tr}&(h(M-N))\\ &\leq a_0n^2 + \sum_{k=1}^{\infty} \frac{2[\mathrm{tr}(M^2) + \mathrm{tr}(D)^2]a_{2k}}{(2k)!}\lambda_\max(L)^{2k-2}\\
    &= a_0n^2 + 2\left(\sum_{k=1}^\infty \frac{a_{2k}}{(2k)!}\lambda_{\max}(L)^{2k-2} \right)[\mathrm{tr}(M^2) + \mathrm{tr}(D)^2]\\ 
    &= \mathrm{tr}\left[a_0 I + 2\left(\sum_{k=1}^\infty \frac{a_{2k}}{(2k)!}\lambda_{\max}(L)^{2k-2} \right) M^2 \right] + C\\
    &\eqqcolon \mathrm{tr}(\tilde h(M)),
\end{align*}
where
\[
C = 2\mathrm{tr}(D)^2\left(\sum_{k=1}^\infty \frac{a_{2k}}{(2k)!}\lambda_{\max}(L)^{2k-2} \right).\qedsymbol
\]
\bigskip
\noindent
\textbf{Proof to Proposition 4.} Since $W_{e_r}$ is doubly stochastic, one has that $W_{e_r}\mathbbm{1} = \mathbbm{1}$ and hence, in the $k$-th iteration,
\begin{align*}
z^{(k, r+1)} &= s^{(k,r+1)} - \bar d\mathbbm{1} = W_{e_r} s^{(k,r)} - W_{e_r} (\bar d\mathbbm{1})\\ &= W_{e_r} (s^{(k,r)} - \bar d\mathbbm{1}) = W_{e_r} z^{(k,r)},
\end{align*}
where 
\[
s^{(k,r)} = \left(\prod_{p=0}^{r-1}W_{e_p}\right)s^{(k)},\quad z^{(k,r)} = s^{(k,r)}-\bar d\mathbbm{1}.
\]
Notice further that $W_{e_{r}}$ is an orthogonal projector $\implies$ $W_{e_r}$ is idempotent and hence
\[
\left\Vert z^{(k, r+1)} \right\Vert_2^2 = z^{(k, r)\top}W_{e_r}z^{(k,r)}.
\]
Taking expectation over all possible edges, one has that
\begin{align*}
\mathbb{E}\left\Vert z^{(k, r+1)} \right\Vert_2^2 &=  z^{(k, r)\top}(\mathbb EW_{e_r})z^{(k,r)}\\ &= z^{(k, r)\top}\left(
I - \frac{1}{2|E|}L
\right)z^{(k,r)}\\ &= \left\Vert z^{(k, r)} \right\Vert_2^2 - \frac{1}{2|E|}z^{(k, r)\top}Lz^{(k,r)}.
\end{align*}
Since $z^{(k,r)}\perp\mathbbm{1}$ by definition, one has that
\[
 \left\Vert z^{(k, r)} \right\Vert_2^2 - \frac{1}{2|E|}z^{(k, r)\top}Lz^{(k,r)} \leq \left(1 - \frac{\lambda_2(L)}{2|E|} \right)\left\Vert z^{(k, r)} \right\Vert_2^2.
\]
And one may then deduce that
\[
\mathbb{E}\left\Vert z^{(k+1)} \right\Vert_2^2 \leq \left(1 - \frac{\lambda_2(L)}{2|E|} \right)^R\left\Vert z^{(k)} \right\Vert_2^2.\qedsymbol
\]

\end{document}